\documentclass[12pt]{article}
\usepackage[fleqn]{amsmath}
\usepackage{amsfonts}
\usepackage{multirow}
\usepackage{verbatim}

\newtheorem{theorem}{Theorem}[section]

\newtheorem{proposition}[theorem]{Proposition}

\DeclareMathOperator{\xNS}{NS}
\DeclareMathOperator{\xPic}{Pic}
\DeclareMathOperator{\xch}{c_1}

\DeclareMathOperator{\Kahler}{K\ddot{a}hler}
\DeclareMathOperator{\xPer}{Per}
\DeclareMathOperator{\xGL}{GL}

\begin{document}

\setlength{\pdfpagewidth}{8.5in}
\setlength{\pdfpageheight}{11in}

\title{Distinguished Line Bundles for Complex Surface Automorphisms}
\author{Paul Reschke}
\date{}
\maketitle

\begin{abstract}
We equate dynamical properties (e.g., positive entropy, existence of a periodic curve) of complex projective surface automorphisms with properties of the pull-back actions of such automorphisms on line bundles. We use the properties of the cohomological actions to describe the measures of maximal entropy for automorphisms with positive entropy.
\end{abstract}

\section{Overview}

In this paper, we investigate cohomological actions induced by complex projective surface automorphisms. In particular, we show that the condition that such an automorphism has positive topological entropy has important implications for the action of the automorphism on line bundles. These constraints lead to a distinguished means of obtaining the measures of maximal entropy for certain projective surface automorphisms. Overall, this paper suggests an initial framework for understanding in greater generality the relationship between cohomological actions and dynamics on K\"{a}hler manifolds. Throughout this paper, the entropy of a map is the topological entropy unless otherwise specified.
\\
\\
A polynomial \(S(t) \in \mathbb{Z}[t]\) is reciprocal if it satisfies \(S(t) = t^s S(t^{-1})\), where \(s\) is the degree of \(S(t)\). A monic irreducible polynomial in \(\mathbb{Z}[t]\) is a Salem polynomial if it is reciprocal, it has a positive real root, and it has exactly two roots with magnitude not equal to one. A real algebraic integer is a Salem number if it is greater than one and its minimal polynomial is a Salem polynomial; the degree of a Salem number is the degree of its minimal polynomial. Given any \(\mathbb{Z}\)-module endomorphism \(\phi\) and any polynomial
\[Q(t) = q_0 + q_1t + \dots + q_nt^n\]
with integer coefficients, there is a naturally defined \(\mathbb{Z}\)-module endomorphism
\[Q(\phi) = q_0\mathbf{1} + q_1\phi + \dots + q_n\phi^n,\]
where \(\mathbf{1}\) is the identity map and \(\phi^j=\phi^{\circ j}\) is the \(j\)-fold iteration \(\phi \circ \cdots \circ \phi\) for any \(j \in \mathbb{N}\).
\begin{theorem}
Let \(X\) be a connected compact K\"{a}hler surface, and let \(\sigma\) be an automorphism of \(X\). Let \(\lambda\) be a Salem number of degree \(s\), let \(S(t)\) be the minimal polynomial for \(\lambda\), and let \(W_S(\sigma)\) be the kernel of the action of \(S(\sigma^*)\) on \(\xPic(X)\). Then the following four statements are equivalent:
\begin{description}
\item[1)] \(W_S(\sigma)\) contains a line bundle with a non-trivial Chern class;
\item[2)] \(W_S(\sigma)\) contains an \(s\)-dimensional sublattice of line bundles with non-trivial Chern classes;
\item[3)] \(X\) is projective and the entropy of \(\sigma\) is \(\log(\lambda) > 0\); and
\item[4)] \(W_S(\sigma)\) contains a nef and big line bundle.
\end{description}
\end{theorem}
Since the entropy of any automorphism of a connected compact \(\Kahler\) surface must be zero or the logarithm of a Salem number, Theorem 1.1 accounts for all smooth complex projective surface automorphisms with positive entropy. Also, any monic irreducible polynomial satisfying (1) in Theorem 1.1 must be either a Salem polynomial or a cyclotomic polynomial. (See \S 2.4 in this paper.)
\\
\\
\textbf{Definition.} Let \(\sigma\) be an automorphism of a smooth complex projective surface \(X\) with entropy \(\log(\lambda)>0\), and let \(S(t)\) be the minimal polynomial for \(\lambda\). Then a line bundle \(L \in \xPic(X)\) will be called \emph{distinguished} if \(S(\sigma^*)L\) is trivial.
\\
\\
By Theorem 1.1, some distinguished line bundle for a smooth complex projective surface automorphism with positive entropy can always be taken to be nef and big; however, there may be no ample distinguished line bundle for such an automorphism.
\begin{theorem}
Let \(X\) be a smooth complex projective surface, and let \(\sigma\) be an automorphism of \(X\) with entropy \(\log(\lambda) > 0\). Then the following two statements are equivalent:
\begin{description}
\item[1)] There is an ample distinguished line bundle on \(X\); and
\item[2)] No curve on \(X\) is periodic for \(\sigma\).
\end{description}
\end{theorem}
Theorem 1.2 shows that, in general, the dichotomy between automorphisms with periodic curves and those without periodic curves can be interpreted as a dichotomy of the sets of distinguished line bundles associated to these automorphisms. A periodic curve may be exceptional (of the first kind), in which case the automorphism will descend to an automorphism of the surface obtained by contraction of the orbit of the exceptional curve. (See \S 2.4 in this paper.) However, an automorphism which is minimal in the sense that it has no periodic exceptional curve may still have a periodic curve. (See \S 6.1, \S 6.2, and \S 6.3 in this paper.)
\\
\\
Many examples of rational surface automorphisms with positive entropy have been shown to exist in \cite{be1}, \cite{dil}, \cite{mc2}, and \cite{ueh}. The constructions in all of these cases yield automorphisms with invariant cuspidal cubic curves; however, examples of rational surface automorphisms with positive entropy and no periodic curves have been show to exist in \(\cite{be2}\). (See \S 6.3 in this paper.) Examples of projective K3 surface automorphisms with positive entropy have been studied in \(\cite{bar}\), \(\cite{ka2}\), \(\cite{lee}\), and \(\cite{mc3}\), among others. All of the examples that arise from torus automorphisms via the Kummer construction must have periodic curves; however, many of the non-Kummer examples have no periodic curves. On the other hand, an abelian surface automorphism with positive entropy can never have a periodic curve. (See \S 6.1 and \S 6.2 in this paper.) As explained in \cite{ca1}, any smooth complex projective surface that admits an automorphism with positive entropy must be birational to an abelian surface, a K3 surface, an Enriques surface, or the projective plane. (See \S 2.4 in this paper.)
\\
\\
Suppose that \(X\) is a smooth complex projective surface, and that \(\sigma\) is an automorphism of \(X\) with positive entropy. Results in \cite{ca2} show that \(\sigma\) has a unique measure of maximal entropy, \(\mu_\sigma\), and that this measure can be expressed as a wedge product \(T_+ \wedge T_-\) of positive closed \((1,1)\)-currents \(T_+\) and \(T_-\) on \(X\) that are, respectively, dilated and contracted under the action on currents induced by \(\sigma\). (See \S 4.1 and \S 4.2 in this paper.)
\begin{theorem}
Let \(X\) be a smooth complex projective surface, and let \(\sigma\) be an automorphism of \(X\) with entropy \(\log(\lambda) > 0\). Suppose that either \(X\) is not a rational surface or no curve on \(X\) is periodic for \(\sigma\). Then there is a distinguished nef and big line bundle on \(X\) whose Chern class contains a semi-positive curvature form; moreover, if \(\omega_0\) is the semi-positive form, then the inductively defined sequence
\[\{\omega_n = (\lambda + \lambda^{-1})^{-1} (\sigma^*\omega_{n-1} + (\sigma^{-1})^*\omega_{n-1})\}_{n \in \mathbb{N}}\]
converges weakly to a positive current \(T\) with the property that the measure \(T \wedge T\) is \(c\mu_\sigma\) for some positive real number \(c\).
\end{theorem}
If \(\sigma\) has no periodic curves in Theorem 1.3, then \(\omega_0\) can be taken to be a \(\Kahler\) form. The assumption that the surface is not rational in Theorem 1.3 allows the use of a theorem due to Kawamata to show the existence of the distinguished semi-positive forms; it is a technical condition at present, as we do not know if there is a rational surface automorphism with positive entropy and no distinguished semi-positive form. (See \S 4.3 in this paper.)
\\
\\
\textbf{Remark on the quadratic case.} Suppose that \(X\) is a smooth complex projective surface, and that \(\sigma\) is an automorphism of \(X\) with no periodic curves whose entropy is the logarithm of a quadratic Salem number. Then \((X;\sigma,\sigma^{-1})\) is a polarized dynamical system of two morphisms in the sense of \cite{ka1} and \cite{lee}; Theorems 1.1 and 1.2 show that this situation is in fact the only one in which a complex surface automorphism gives rise to such a dynamical system. In this special setting, results from \cite{ka1}, \cite{lee}, and \cite{ca2} guarantee the existence of the \(\Kahler\) form in Theorem 1.3; thus Theorem 1.3 is a broad generalization of the previously developed descriptions of measures of maximal entropy for polarized dynamical systems of two automorphisms. Related observations about the quadratic case appear in \cite{ka2} and \cite{cal}.
\\
\\
The methods used to characterize distinguished line bundles for Theorem 1.1 also provide a means of understanding in general the orbits of line bundles under complex projective surface automorphisms with positive entropy.
\begin{theorem}
Let \(X\) be a smooth complex projective surface, and let \(\sigma\) be an automorphism of \(X\) with entropy \(\log(\lambda) > 0\). Then the following two statements are equivalent:
\begin{description}
\item[1)] There is an ample line bundle on \(X\) whose orbit of Chern classes under \(\sigma^*\) spans \(\xNS(X) \otimes \mathbb{R}\); and
\item[2)] The characteristic polynomial for the action of \(\sigma^*\) on \(\xNS(X)\) is separable.
\end{description}
\end{theorem}
In the setting of Theorem 1.4, the characteristic polynomial for the action of \(\sigma^*\) on \(\xNS(X)\) may be separable or inseparable; there are examples for both cases. (See \S 6.1, \S 6.2, and \S 6.3 in this paper.) 
\\
\\
The remainder of this paper is organized as follows: in \S 2, we describe important properties of cohomological actions induced by automorphisms of compact \(\Kahler\) surfaces; in \S 3, we use these properties to prove Theorems 1.1 and 1.2; in \S 4, we discuss measures of maximal entropy and we prove Theorem 1.3; in \S 5, we discuss cyclotomic factors for surface automorphisms and we prove Theorem 1.4; finally, in \S 6, we give examples of automorphisms that highlight the different situations that can arise with regard to induced actions on line bundles.
\\
\\
\emph{Acknowledgements:} We thank Laura DeMarco, Serge Cantat, Izzet Coskun, and Chong Gyu Lee for helpful discussions. We also thank Igor Dolgachev for asking a question that prompted the development of Theorem 1.4 and the material in \S 5. 

\section{Cohomological Structures}

In this section, we recall a variety of facts about cohomology groups (\S 2.1), line bundles (\S 2.2), and differential forms (\S 2.3) on compact \(\Kahler\) surfaces. (For detailed background material, see \cite{bhpv}, \cite{har}, and \cite{huy}.) In \S 2.4, we review several previously known results about the implications of the existence of an automorphism with positive entropy on a compact \(\Kahler\) surface.
\\
\\
Let \(X\) be a connected compact \(\Kahler\) surface, and let \(\sigma\) be an automorphism of \(X\).

\subsection{Cohomology Groups}

The Picard group of \(X\) is denoted \(\xPic(X)\), and the N\'eron-Severi group of \(X\) is denoted \(\xNS(X)\); the first Chern map from \(\xPic(X)\) to \(H^2(X,\mathbb{Z})\) is denoted \(\xch\). When torsion is factored out, the image of \(\xch\) is \(\xNS(X)\); the Lefschetz theorem on (1,1) classes gives
\[\xNS(X) = i(H^2(X,\mathbb{Z})) \cap H^{1,1}(X),\]
where \(i\) is the natural map from \(H^2(X,\mathbb{Z})\) to \(H^2(X,\mathbb{C})\). The kernel of \(\xch\) (i.e., the Picard variety of \(X\)) is denoted \(\xPic^0(X)\).
\\
\\
The Hodge decomposition gives
\[H^2(X,\mathbb{R}) = H^2(X,\mathbb{Z}) \otimes \mathbb{R} = H^2(X,\mathbb{C})_\mathbb{R} = H^{1,1}(X)_\mathbb{R} \oplus (H^{2,0}(X) \oplus H^{0,2}(X))_\mathbb{R}.\]
The cup product defines a quadratic form on \(H^2(X,\mathbb{R})\); the restriction of this form to \(H^2(X,\mathbb{Z})\) coincides with the image of the intersection form on \(\xPic(X)\). The pull-back map \(\sigma^*\) on each of these spaces is an automorphism that preserves the pairing; moreover, \(\sigma^*\) commutes with \(\xch\). The intersection of two elements \(o_1\) and \(o_2\) is denoted \(o_1.o_2\), while the self-intersection of an element \(o\) is denoted \(o^2\).

\subsection{Line Bundles}

A line bundle \(L \in \xPic(X)\) is ample if and only if \(L^2 > 0\) and \(L.[D] > 0\) for any effective divisor \(D\) on \(X\); \(X\) is projective if and only if there is a line bundle \(L \in \xPic(X)\) with \(L^2 > 0\). (See \cite{bhpv}, \S IV.6.) If \(X\) is projective, then \(\xPic(X)\) is precisely the group of divisors on \(X\) modulo linear equivalence. (See \cite{har}, \S II.4 and \S II.6.) The set of effective divisor classes on \(X\) is preserved by \(\sigma^*\), as is the set of ample line bundles on \(X\). (See \cite{bhpv}, \S I.6, and \cite{har}, \S II.6.) The following property of ample line bundles is a consequence of the Hodge index theorem: if \(L \in \xPic(X)\) is ample and \(L' \in \xPic(X)\) satisfies \(L.L' = 0\), then \((L')^2 \leq 0\). (See \cite{har}, \S V.1.)
\\
\\
A line bundle \(L \in \xPic(X)\) is called nef if \(L.[D] \geq 0\) for any effective divisor \(D\) on \(X\); if \(X\) is projective, a nef line bundle on \(X\) is called big if it has positive self-intersection. (See \cite{bhpv}, \S I.6, \S IV.7, and \S IV.12.) Thus any ample line bundle on \(X\) is necessarily nef and big. The following property of line bundles with positive self-intersection is a consequence of the Riemann-Roch theorem: if \(L \in \xPic(X)\) has \(L^2 > 0\) and \(L.H > 0\) for some ample \(H \in \xPic(X)\), then \(L^{\otimes m}\) is effective for some \(m \in \mathbb{N}\). (See \cite{har}, \S V.1.)

\subsection{Differential Forms}

The space of complex differential forms on \(X\) is
\[\Omega(X) = \bigoplus\nolimits_{0 \leq r \leq 4} \Omega^r(X) = \bigoplus\nolimits_{0 \leq r \leq 4} \left( \bigoplus\nolimits_{p+q=r} \Omega^{p,q}(X) \right),\]
where each \(\Omega^{p,q}(X)\) is the space of complex differential \((p,q)\)-forms on \(X\). A complex differential form is real if it is equal to its own complex conjugate.
\\
\\
The signature of the quadratic form on \(H^{1,1}(X)_\mathbb{R}\) is \((1,h^{1,1}-1)\). This fact leads to an analogue of the Hodge index theorem that applies to real closed \((1,1)\)-forms: if \(v_1\) and \(v_2\) are linearly independent in \(H^{1,1}(X)_\mathbb{R}\) with \(v_1^2 > 0\) and \(v_2^2 \geq 0\), then \(v_1.v_2 \neq 0\).
\\
\\
A real \((1,1)\)-form \(\omega\) on \(X\) is (semi-)positive if it prescribes a positive (semi-)definite Hermitian form on the tangent space at every point in \(X\). A \(\Kahler\) form is a real \((1,1)\)-form that is positive and \(d\)-closed; a \(\Kahler\) form must have positive intersection with every effective divisor on \(X\). The set of \(\Kahler\) classes forms a convex cone \(C_K(X) \subseteq H^{1,1}(X)_\mathbb{R}\); it is contained in \(C_+(X) \subseteq H^{1,1}(X)_\mathbb{R}\), the positive cone of \(X\), which is also convex. The \(\Kahler\) cone is open in \(H^{1,1}(X)_\mathbb{R}\), and, it is given explicitly by
\[C_K(X) =\{v \in C_+(X) \mid v.\xch([Y])>0 \, \forall \, \text{curve} \, Y \subseteq X \, \text{with} \, Y^2 < 0\}.\]
(See \cite{buc}, \cite{dem}, \cite{lam}, and \cite{voi}, \S 1.) Thus the Chern class of a line bundle \(L \in \xPic(X)\) is represented by a \(\Kahler\) form if and only if \(L\) is ample. (This is the content of the Kodaira embedding thoerem; see \cite{bal}, \S 9, and \cite{voi}, \S 1.) So \(X\) is projective if and only if \(C_K(X) \cap \xNS(X) \neq \emptyset\). Both \(C_K(X)\) and \(C_+(X)\) are preserved by \(\sigma^*\).

\subsection{Positive Entropy}

The topological entropy of \(\sigma\) is \(\log(\lambda)\), where \(\lambda\) is the spectral radius of the linear map
\[\sigma^*:H^{1,1}(X)_\mathbb{R} \rightarrow H^{1,1}(X)_\mathbb{R}.\]
(See \cite{ca1}, \S 4.4.2.) A real algebraic integer is called a Salem number if it is greater than one and it has a reciprocal minimal polynomial with exactly two roots off the unit circle. It is a well-known fact that the entropy of \(\sigma\) is either zero or the logarithm of a Salem number; indeed, the irreducible factors of the characteristic polynomial for the action of \(\sigma^*\) on \(H^2(X,\mathbb{Z})\) can only consist of cyclotomic polynomials and at most one Salem polynomial (counting multiplicity). (See, e.g., \cite{mc1}, \cite{res}, and \cite{ueh}.) The degree of a Salem number is the degree of its minimal polynomial; it is necessarily even. If the entropy of \(\sigma\) is \(\log(\lambda)\), then (for all \(k \in \mathbb{N}\)) the entropy of \(\sigma^k\) is \(\log(\lambda^k)\). If \(\lambda\) is a Salem number of degree \(s\), then \(\lambda^k\) is a Salem number of degree \(s\) as well; indeed, the set of Galois conjugates of \(\lambda^k\) is precisely the set of all \(k\)-th powers of Galois conjugates of \(\lambda\).
\\
\\
The condition that \(\sigma\) has positive entropy constrains what \(X\) can be: after any \((-1)\)-curves with finite order under \(\sigma\) are contracted, \(X\) must be a complex torus, a K3 surface, an Enriques surface, or a rational surface; moreover, if \(X\) is a rational surface, then \(X\) must be a blow-up of the projective plane at ten or more points. (See \cite{ca1}, \S 10.3.) If there is a \((-1)\)-curve \(E\) on \(X\) that is periodic for \(\sigma\), then \(\sigma\) descends to an automorphism \(\sigma'\) of the surface \(X'\) obtained from \(X\) by contraction of the orbit of \(E\); moreover, since the orbit of \(\xch([E])\) in \(\xNS(X) \subseteq H^{1,1}(X)_\mathbb{R}\) contributes nothing to the entropy of \(\sigma\), \(\sigma'\) and \(\sigma\) must have the same entropy. (See also \cite{ka2}, \S 2, and \cite{ca1}, \S 4.1.) Conversely, if \(x \in X\) is a periodic point for \(\sigma\), then \(\sigma\) extends (without a change in entropy) to an automorphism of the blow-up of \(X\) at the points in the orbit of \(x\).

\section{Distinguished Line Bundles}

In this section, we prove Theorems 1.1 and 1.2; the proofs are given in \S 3.5, and they directly invoke the propositions in the preceding subsections.

\subsection{From Distinguished Chern Classes to Positive Entropy}

\begin{proposition}
Let \(X\) be a connected compact K\"{a}hler surface, and let \(\sigma\) be an automorphism of \(X\). Let \(\lambda\) be a Salem number, and let \(S(t)\) be the minimal polynomial for \(\lambda\). Suppose that there is a non-trivial Chern class \(\xch(L) \in \xNS(X)\) such that \(S(\sigma^*)\xch(L)=0\). Then \(X\) is projective and \(\sigma\) has entropy \(\log(\lambda)\).
\end{proposition}

\emph{Proof:} The set of eigenvalues of the linear action of \(S(\sigma^*)\) on \(\xNS(X) \otimes \mathbb{R}\) is precisely the set of all \(S(\alpha)\) where \(\alpha\) is an eigenvalue of the action of \(\sigma^*\). Since zero is an eigenvalue of \(S(\sigma^*)\), the eigenvalues of \(\sigma^*\) must include a root of \(S(t)\); moreover, since the characteristic polynomial for the action of \(\sigma^*\) on \(\xNS(X) \otimes \mathbb{R}\) has integer coefficients, it must have \(S(t)\) as a factor. So the entropy of \(\sigma\) must be \(\log(\lambda)\). Let \(D_+\) be the eigenspace for \(\sigma^*\) corresponding to the eigenvalue \(\lambda\), let \(D_-\) be the eigenspace for \(\sigma^*\) corresponding to the eigenvalue \(\lambda^{-1}\), and let \(E = D_+ \oplus D_-\). Then, since the dimension of \(E\) is two and the signature of \(H^{1,1}(X)_\mathbb{R}\) is \((1,h^{1,1}-1)\), \(E\) cannot be a totally isotropic subspace of \(\xNS(X) \otimes \mathbb{R}\). Since any vector in \(D_+\) or \(D_-\) must have zero self-intersection, it follows that the signature of \(E\) must be \((1,1)\). If the self-intersection of every element in \(\xNS(X)\) were non-positive, then the same would be true for every element in \(\xNS(X) \otimes \mathbb{Q}\), which is dense in \(\xNS(X) \otimes \mathbb{R}\); but then, by continuity, the self-intersection of every element in \(\xNS(X) \otimes \mathbb{R}\) would be non-positive, which is not the case. Thus \(\xNS(X)\) must contain some Chern class \(\xch(L_0)\) with
\[\xch(L_0)^2 = L_0^2 > 0,\] and \(X\) must be projective. \(\Box\)

\subsection{From Positive Entropy to Distinguished Chern Classes}

\begin{proposition}
Let \(X\) be a connected compact K\"{a}hler surface, and let \(\sigma\) be an automorphism of \(X\). Let \(\lambda\) be a Salem number, let \(S(t)\) be the minimal polynomial for \(\lambda\), and let \(s\) be the degree of \(S(t)\). Suppose that \(X\) is projective and that the entropy of \(\sigma\) is \(\log(\lambda)\). Then there is an \(s\)-dimensional sublattice of \(\xNS(X)\) that is annihilated by \(S(\sigma^*)\).
\end{proposition}

\emph{Proof:} Let \(n \in \mathbb{N}\) be the dimension of \(H^{1,1}(X)_\mathbb{R}\). So the eigenvalues of \(\sigma^*\) acting on \(H^{1,1}(X)_\mathbb{R}\) are \(\lambda\), \(\lambda^{-1}\), and \(n-2\) algebraic integers with magnitude one (counting multiplicity). Let \(D_+\) be the eigenspace corresponding to \(\lambda\), let \(D_-\) be the eigenspace corresponding to \(\lambda^{-1}\), and let \(E = D_+ \oplus D_-\). So \(E\) has signature \((1,1)\), \(E^\perp\) has signature \((0,n-2)\), and \(H^{1,1}(X)_\mathbb{R} = E \oplus E^\perp\). If \(\lambda\) were not an eigenvalue of \(\sigma^*\) acting on \(\xNS(X)\), then \(\xNS(X) \otimes \mathbb{R}\) would necessarily be contained in \(E^\perp\); but then \(\xNS(X) \otimes \mathbb{R}\) would be negative definite and \(X\) could not be projective. So \(\xNS(X) \otimes \mathbb{R}\) must contain both \(D_+\) and \(D_-\), and \(S(t)\) must be a factor in the characteristic polynomial for \(\sigma^*\) acting on \(\xNS(X) \otimes \mathbb{R}\).
\\
\\
{\bf Lemma 3.2.1} \emph{Let \(\phi\) be an invertible linear transformation of a vector space \(V\) over \(\mathbb{Q}\), and suppose that \(R(t)\) is an irreducible factor (over \(\mathbb{Q}\)) of the characteristic polynomial for \(\phi\) while \(R(t)^2\) is not a factor. Then there is a subspace \(V' \subseteq V\) (defined over \(\mathbb{Q}\)) such that \(\phi(V')=V'\) and the characteristic polynomial for \(\phi |_{V'}\) is \(R(t)\).}
\\
\\
\emph{Proof of Lemma 3.2.1:} Consider \(V\) as an \(\mathbb{Q}[t]\)-module where the action of \(t\) on \(V\) is the application of \(\phi\); then the decomposition of \(V\) with respect to elementary divisors is preserved by \(\phi\) and has \(\mathbb{Q}[t]/(R(t))\) as a factor--and the characteristic polynomial for \(\phi\) on \(\mathbb{Q}[t]/(R(t))\) is \(R(t)\). (This process is analogous to that of expressing \(\phi\) in rational canonical form--where \(V\) is decomposed according to invariant factors rather than elementary divisors; see, e.g., \cite{dum}, \S 12.) \(\Box\)
\\
\\
With \(\phi=\sigma^*\) and \(V = \xNS(X) \otimes \mathbb{Q}\), Lemma 3.2.1 gives an \(s\)-dimensional subspace \(V' \subseteq \xNS(X) \otimes \mathbb{Q}\) that is annihilated by \(S(\sigma^*)\). Thus, since any element of \(\xNS(X) \otimes \mathbb{Q}\) has some multiple which is an element of \(\xNS(X)\), \(V' \cap \xNS(X)\) is an \(s\)-dimensional sublattice of \(\xNS(X)\) that is annihilated by \(S(\sigma^*)\). \(\Box\)

\subsection{From Chern Classes to Line Bundles}

Let \(X\) be a connected compact \(\Kahler\) surface, and let \(\sigma\) be an automorphism of \(X\). If \(X\) is a finite blow-up of a K3 surface, an Enriques surface, or a rational surface, then the first Betti number of \(X\) is zero; so the dimension of \(H^1(X,\mathcal{O}_X)\) is zero, and
\[\xch: \xPic(X) \rightarrow \xNS(X)\]
is an isomorphism of non-torsion elements. (See \cite{bhpv}, \S I.9, \S III.4, \S V.1, \S VIII.2, and \S VIII.15.) If \(X\) is not one of these three types of surfaces and \(\sigma\) has positive entropy, then \(X\) must be a finite blow-up of a torus; so the dimension of \(H^1(X,\mathcal{O}_X)\) is two. In any case, for any polynomial \(S(t) \in \mathbb{Z}[t]\), if a line bundle \(L \in \xPic(X)\) satisfies \(S(\sigma^*)L=0\), then the Chern class \(\xch(L)\) also satisfies \(S(\sigma^*)\xch(L)=0\). 
\begin{proposition}
Let \(X\) be a connected compact K\"{a}hler surface, and let \(\sigma\) be an automorphism of \(X\). Let \(\lambda\) be a Salem number, and let \(S(t)\) be its minimal polynomial. Suppose that the entropy of \(\sigma\) is \(\log(\lambda)\) and that \(L_2 \in \xPic(X)\) is a line bundle satisfying \(S(\sigma^*)\xch(L_2)=0\). Then there is a line bundle \(L_1 \in \xPic(X)\) that satisfies \(S(\sigma^*)L_1=0\) and \(\xch(L_1) = \xch(L_2)\).
\end{proposition}

\emph{Proof:} If \(X\) is not bimeromorphic to a torus, then the statement is evident, and, moreover, \(L_1=L_2\).

If \(X\) is a torus, then \(H^*(X,\mathbb{Z})\) is generated by \(H^1(X,\mathbb{Z})\) via the cup product; so, in particular, the six eigenvalues for the action of \(\sigma^*\) on \(H^2(X,\mathbb{Z})\) are precisely the products of all pairs among the four eigenvalues for the action of \(\sigma^*\) on \(H^1(X,\mathbb{Z})\). (See \cite{res}, \S 2.) Since the eigenvalues for the action of \(\sigma^*\) on \(H^2(X,\mathbb{Z})\) are \(\lambda\), \(\lambda^{-1}\), and four algebraic integers with magnitude one, the eigenvalues for the action of \(\sigma^*\) on \(H^1(X,\mathbb{Z})\) must be two algebraic integers with magnitude \(\sqrt{\lambda}\) and two algebraic integers with magnitude \(\sqrt{\lambda}^{-1}\); so, since no root of \(S(t)\) has magnitude \(\sqrt{\lambda}\) or \(\sqrt{\lambda}^{-1}\), the action of \(S(\sigma^*)\) on \(H^1(X,\mathbb{R})\) must be surjective. Moreover, since \(\xPic^0(X)\) is the quotient of \(H^1(X,\mathcal{O}_X)\) by an embedding of \(H^1(X,\mathbb{Z})\) (that commutes with \(\sigma^*\)), the action of \(S(\sigma^*)\) must in fact be surjective on \(\xPic^0(X)\). Thus there is some \(L_0 \in \xPic^0(X)\) such that
\[S(\sigma^*)L_0 = S(\sigma^*)L_2,\]
and the line bundle \(L_1 = L_2 - L_0\) satisfies \(S(\sigma^*)L_1=0\) with \(\xch(L_1) = \xch(L_2)\). 

It remains true that \(S(\sigma^*)\) must be surjective on \(\xPic^0(X)\) if \(X\) is a finite blow-up of a torus, so that \(L_0\) and \(L_1\) can be constructed as above in this case as well. \(\Box\)

\subsection{Nef and Big Line Bundles}

Let \(X\) be a smooth complex projective surface, and suppose that \(\sigma\) is an automorphism of \(X\) with entropy \(\log(\lambda)>0\). Let \(S(t)\) be the minimal polynomial for \(\lambda\), and let \(s\) be the degree of \(S(t)\). For the action of \(\sigma^*\) on \(H^{1,1}(X)_\mathbb{R}\), let \(D_+\) be the eigenspace corresponding to \(\lambda\), let \(D_-\) be the eigenspace corresponding to \(\lambda^{-1}\), and let \(E = D_+ \oplus D_-\). Let \(n \in \mathbb{N}\) be the dimension of \(H^{1,1}(X)_\mathbb{R}\), and let \(\{v_1,\dots,v_n\}\) be a basis for \(H^{1,1}(X)_\mathbb{R}\) such that \(v_1 \in D_+\), \(v_2 \in D_-\), and \(v_j \in E^\perp\) for all other \(j\). Since \(C_K(X)\) is open in \(H^{1,1}(X)_\mathbb{R}\), it must contain some element \(v\) whose first two coordinates are some non-trivial elements \(e_+ \in D_+\) and \(e_- \in D_-\). Then
\[\lim_{k \to \infty} \lambda^{-k}(\sigma^*)^k v = e_+\]
and
\[\lim_{k \to \infty} \lambda^{-k}((\sigma^{-1})^*)^k v = e_-.\]
Thus, since \(\sigma^*\) preserves the \(\Kahler\) cone, \(\overline{C_K(X)}\) must contain \(e_+\) and \(e_-\), as well as \(ae_+ + be_-\) for any two non-negative real numbers \(a\) and \(b\).
\\
\\
Let \(u=ae_+ + be_-\), with \(a\) and \(b\) any two positive real numbers. If \(u^2\) were zero, then \(E\) would be a two-dimensional totally isotropic subspace of \(H^{1,1}(X)_\mathbb{R}\), which cannot be the case; so \(u^2\) is positive. Also, \(u.\xch([D])\) is non-negative for any effective divisor \(D\) on \(X\). (See also \cite{ka2}, \S 2.)
\begin{proposition}
[\cite{ka2}, Proposition 3.1] Assume the hypotheses and notation of the preceding text in this subsection, and let \(C\) be an irreducible curve on \(X\). Then \(C\) is periodic for \(\sigma\) if and only if \(u.\xch([C])\) is zero. There are only finitely many such curves on \(X\).
\end{proposition}
Thus, in particular, \(u\) is a \(\Kahler\) class if and only if no curve on \(X\) is periodic for \(\sigma\). If \(C\) is a curve on \(X\) with \(u.\xch([C])=0\), then, by the analogue of the Hodge index theorem, \(C\) must have negative self-intersection.
\\
\\
Let \(E' \subseteq \xNS(X) \otimes \mathbb{Q}\) be the \(s\)-dimensional subspace that is annihilated by \(S(\sigma^*)\); so \(E\) is contained in \(E' \otimes \mathbb{R}\). Let \(\xNS'(X)\) be the \(s\)-dimensional sublattice \(E' \cap \xNS(X)\); so \(\xNS'(X)\) is precisely the set of all Chern classes in \(\xNS(X)\) that are annihilated by \(S(\sigma^*)\). Since \(E^\perp\) is negative definite, \(\xNS'(X) \cap C_+(X)\) must contain some non-trivial Chern class \(\xch(L_+)\). In light of Proposition 3.3, take \(L_+\) to satisfy \(S(\sigma^*)L_+=0\).
\\
\\
Suppose that \(C\) is an irreducible curve on \(X\) such that \(\sigma^k(C)=C\) for some \(k \in \mathbb{N}\), and let \(S_k(t)\) be the minimal polynomial for \(\lambda^k\); so \(S_k(t)\) is again a Salem polynomial of degree \(s\). Since the entropy of \(\sigma^k\) is \(\log(\lambda^k)\) and every invariant space for \(\sigma^*\) is also invariant for \((\sigma^*)^k=(\sigma^k)^*\), it follows that \(\xNS'(X)\) is annihilated by \(S_k((\sigma^*)^k)\). Thus, for any \(\xch(L) \in \xNS'(X)\),
\[0 = (S_k((\sigma^*)^k)[C]).(S_k((\sigma^*)^k)L) = (\Sigma^2)([C].L),\]
where \(\Sigma\) is the sum of the coefficients of \(S_k(t)\); but since one is not a root of \(S_k(t)\), \(\Sigma\) cannot be zero--which forces \([C].L = 0\). So, in particular, \(\xNS'(X) \cap C_K(X)\) must be empty if \(\sigma\) has any periodic curves.
\begin{proposition}
Assume the hypotheses and notation of the preceding text in this subsection, and suppose that no curve on \(X\) is periodic for \(\sigma\). Then there is an ample line bundle \(L \in \xPic(X)\) satisfying \(S(\sigma^*)L=0\).
\end{proposition}

\emph{Proof:} Let \(\{w_3,\dots,w_s\}\) be a basis for \(E^\perp\) in \(E' \otimes \mathbb{R}\); so \(\{v_1,v_2,w_3,\dots,w_s\}\) is a basis for \(E' \otimes \mathbb{R}\). Since \((L_+)^2\) is positive, the first two coordinates of \(\xch(L_+)\) must be non-zero; moreover, since \(\xch(L_+)\) has non-negative intersection with both \(e_+\) and \(e_-\), the first two coordinates of \(\xch(L_+)\) must be \(ae_+\) and \(be_-\) for some positive numbers \(a\) and \(b\). So
\[\lim_{k \to \infty} \lambda^{-k}((\sigma^*)^k+((\sigma^{-1})^*)^k) \xch(L_+) = ae_+ + be_-,\]
which is a \(\Kahler\) class because \(\sigma\) has no periodic curves. Since \(C_K(X)\) is open, there is some positive integer \(k'\) such that the \(k'\)-th iterate of the sequence is also a \(\Kahler\) class. So \(((\sigma^*)^{k'}+((\sigma^{-1})^*)^{k'}) \xch(L_+)\) is a \(\Kahler\) class and an element of \(\xNS'(X)\), and \(((\sigma^*)^{k'}+((\sigma^{-1})^*)^{k'}) L_+\) is an ample line bundle that is annihilated by \(S(\sigma^*)\). \(\Box\)
\\
\\
If no curve on \(X\) is periodic for \(\sigma\) and the degree of \(S(t)\) is two, then \(\xch(L_+)\) is an element of \(E \cap C_K(X)\) and \(L_+\) itself is ample. In general (whether or not \(\sigma\) has periodic curves), if the degree of \(S(t)\) is two, then \(\xch(L_+)\) is an element of \(E \cap \overline{C_K(X)}\) and \(L_+\) is nef and big. (See also \cite{ka2}, \S 3.)
\\
\\
Suppose that \(C\) is an irreducible curve on \(X\) with non-negative self-intersection; so \(\xch([C])\) is an element of \(\overline{C_+(X)}\). Then, since the intersection of any two elements in \(C_+(X)\) is non-negative, the intersection of \(\xch([C])\) with any element of \(C_+(X)\) must be non-negative. Thus the only barrier to the existence of a nef and big line bundle whose Chern class is contained in \(\xNS'(X)\) is the set of irreducible curves on \(X\) with negative self-intersection.
\begin{proposition}
Assume the hypotheses and notation of the preceding text in this subsection. Then there is a nef and big line bundle \(L \in \xPic(X)\) satisfying \(S(\sigma^*)L=0\).
\end{proposition}

\emph{Proof:} For some \(q \in \mathbb{N}\), \((L_+)^{\otimes q}\) is an effective divisor class; let \(Y_+\) be an effective divisor in this class. Let \(B \in \mathbb{Z}\) be the minimum value of the self-intersection of an irreducible curve in the support of \(Y_+\); so, for any \(k \in \mathbb{N}\), \(B\) is a lower bound for the self-intersection of an irreducible curve in the support of
\[((\sigma^*)^k + ((\sigma^{-1})^*)^k)Y_+.\]
If an irreducible curve has negative intersection with an effective divisor, then the curve must be in the support of the divisor; thus, for any \(k \in \mathbb{N}\),
\[(((\sigma^*)^{k}+((\sigma^{-1})^*)^{k})((L_+)^{\otimes q})\]
has non-negative intersection with every irreducible curve on \(X\) with either non-negative self-intersection or self-intersection less than \(B\). Let \(\mathcal{C}\) be the set of all irreducible curves on \(X\) with self-intersection at least \(B\) and at most \(-1\). Let \(w\) be the limit of the sequence
\[\{\lambda^{-k}((\sigma^*)^k+((\sigma^{-1})^*)^k) \xch((L_+)^{\otimes q})\}_{k \in \mathbb{N}};\]
so \(w=ae_++be_-\) for some positive numbers \(a\) and \(b\), and hence has positive self-intersection and non-negative intersection with every curve on \(X\). Let \(w'=ae_+-be_-\), let \(m \in \mathbb{N}\) be the dimension of \(\xNS(X) \otimes \mathbb{R}\), and let \(\{w_3,\dots,w_m\}\) be a basis for \(E^\perp\) in \(\xNS(X) \otimes \mathbb{R}\); so \(\{w,w'\}\) is a basis for \(E\), \(\{w,w',w_3,\dots,w_m\}\) is a basis for \(\xNS(X) \otimes \mathbb{R}\), and \(\{w',w_3,\dots,w_m\}\) is a basis for \(<w>^\perp\) in \(\xNS(X) \otimes \mathbb{R}\). Thus, since the signature of \(\xNS(X) \otimes \mathbb{R}\) is \((1,m-1)\),  for each \(K \in \{B,\dots,-1\}\) the set
\[\{v \in \xNS(X) \otimes \mathbb{R} \mid v^2 = K, 0 \leq v.w \leq 1\}\]
is homeomorphic to the set
\[\{\vec{x} \in \mathbb{R}^{m-1} \mid -K \leq ||\vec{x}||^2 \leq 1-K\},\]
and hence is compact. So, in particular, the set of all Chern classes of curves in \(\mathcal{C}\) intersecting \(w\) with value one or less is finite. Let \(\mathcal{C}_0\) be the set of curves in \(\mathcal{C}\) whose Chern classes have intersection zero with \(w\); then every curve in \(\mathcal{C}_0\) is periodic, and hence must have intersection zero with every element of \(\xNS'(X)\). Let \(\epsilon>0\) be the minimum value of the intersection of \(w\) with a Chern class of a curve in \(\mathcal{C}-\mathcal{C}_0\). Suppose that \(\{C_l\}_{l \in \mathbb{N}}\) is a sequence of curves in \(\mathcal{C}-\mathcal{C}_0\); then the sequence
\[\{y_l = (1/(\xch([C_l]).w))\xch([C_l])\}_{l \in \mathbb{N}}\]
is contained in the compact set
\[\{v \in \xNS(X) \otimes \mathbb{R} \mid v.w=1, B/\epsilon^2 \leq v^2 \leq 0\},\]
and hence must have a subsequence converging to some element \(y \in \xNS(X) \otimes \mathbb{R}\). If there were also a sequence \(\{x_l\}_{l \in \mathbb{N}}\) of elements in \(\xNS(X) \otimes \mathbb{R}\) converging to \(w\) such that \(\xch([C_l]).x_l \leq \epsilon/2\) for every \(l\), then it would follow that \(y.w = 1\) while \(y_l.x_l \leq 1/2\) for every \(l\)--which cannot happen. So there is an open neighborhood \(U\) of \(w\) in \(\xNS(X) \otimes \mathbb{R}\) such that every element of \(U\) intersects every Chern class of a curve in \(\mathcal{C}-\mathcal{C}_0\) with value greater than \(\epsilon/2\), and thus there is some positive integer \(k'\) such that
\[(((\sigma^*)^{k'}+((\sigma^{-1})^*)^{k'}) L_+).[C] > \frac{\lambda^{k'} \epsilon}{2} \; \; \; \forall \, C \in \mathcal{C}-\mathcal{C}_0.\]
So \(((\sigma^*)^{k'}+((\sigma^{-1})^*)^{k'}) L_+\) is a nef and big line bundle that is annihilated by \(S(\sigma^*)\). \(\Box\)

\subsection{Proofs of Theorems}

\emph{Proof of Theorem 1.1:} By Proposition 3.1, (1) implies (3); by Propositions 3.2 and 3.3, (3) implies (2); and it is evident that (2) implies (1). Proposition 3.6 then shows that any of (1), (2), or (3) imply (4); and it is evident that (4) implies (1). \(\Box\)
\\
\\
\emph{Proof of Theorem 1.2:} Proposition 3.5 is one direction; the other direction follows immediately from the discussion preceding Proposition 3.5. \(\Box\)

\section{Measures of Maximal Entropy}

In this section, we discuss various means of characterizing the unique measure of maximal entropy for a projective surface automorphism with positive entropy. We prove Theorem 1.3 in \S 4.3, drawing on material from \S 4.1; in \S 4.2, we discuss the role that subvarieties play in the understanding of measures of maximal entropy. (For detailed background material on topological entropy and measures of maximal entropy, see \cite{pet}.)

\subsection{Positive Currents}

Suppose that \(X\) is a connected compact \(\Kahler\) surface, and that \(\sigma\) is an automorphism of \(X\); so the topological entropy of \(\sigma\) is equal to the logarithm of the spectral radius of the action on \(H^{1,1}(X)_\mathbb{R}\) induced by \(\sigma\). 
\\
\\
An \(r\)-current on \(X\) is a real-valued linear functional on \(\Omega^{4-r}(X)_\mathbb{R}\) that is continuous with respect to the topology of uniform convergence. The space of currents on \(X\) is
\[W(X) = \bigoplus\nolimits_{0 \leq r \leq 4}W^r(X),\]
where each \(W^r(X)\) is the space of \(r\)-currents on \(X\).  (See also \cite{lel}, \S I.1, and \cite{bhpv}, \S I.11.) The exterior derivative \(d\) on differential forms induces a map on currents: if \(T\) is an \(r\)-current, then \(dT\) is the \(r+1\) current given by, for \(\eta \in \Omega^{3-r}(X)\),
\[dT(\eta) = (-1)^{r+1}T(d\eta).\]
Since \(d^2=0\) as a map on currents, \(d\) gives rise to cohomology groups \(H^r_W(X)\). For any \(\omega \in \Omega^r(X)_\mathbb{R}\), let \(T_\omega\) be the \(r\)-current given by, for \(\eta \in \Omega^{4-r}(X)_\mathbb{R}\), \[T_\omega(\eta) = \int_X \omega \wedge \eta;\]
then the map from \(\Omega(X)_\mathbb{R}\) to \(W(X)\) given by \(\omega \mapsto T_\omega\) commutes with \(d\) and induces isomorphisms from the De Rham cohomology groups \(H^r(X,\mathbb{R})\) to the groups \(H^r_W(X)\). Any \(d\)-closed \(r\)-current on \(X\) descends to a linear functional on \(H^{4-r}(X,\mathbb{R})\); two \(d\)-closed \(r\)-currents represent the same cohomology class in \(H^r_W(X)\) if and only if they give the same linear functional on \(H^{4-r}(X,\mathbb{R})\). For a \(d\)-closed current \(T \in W^2(X)\) and a \(d\)-closed form \(\omega \in \Omega^2(X)_\mathbb{R}\), the intersection pairing on \(H^2(X,\mathbb{R})\) is given by \([T].[\omega] = T(\omega)\). The push-forward by \(\sigma\) of an \(r\)-current \(T\) is the \(r\)-current \(\sigma_*T\) given by, for \(\eta \in \Omega^{4-r}(X)_\mathbb{R}\), \(\sigma_*T(\eta) = T(\sigma^*\eta)\); the pull-back of an \(r\)-current \(T\) by \(\sigma\) is the \(r\)-current \(\sigma^*T = (\sigma^{-1})_*T\). The pull-back map \(\sigma^*\) commutes with the map from \(\Omega(X)_\mathbb{R}\) to \(W(X)\) given by \(\omega \mapsto T_\omega\).
\\
\\
A \((1,1)\)-current on \(X\) is a \(2\)-current which takes non-zero values only on elements of \(\Omega^{1,1}(X)_\mathbb{R}\); every \(2\)-current on \(X\) decomposes uniquely as the sum of a \((1,1)\)-current and a \(2\)-current whose restriction to \(\Omega^{1,1}(X)_\mathbb{R}\) is trivial. (See also \cite{ca1}, \S 5.1, and \cite{ka1}, \S 3.2.) The cohomology class of any \(d\)-closed \((1,1)\)-current is a class in \(H^{1,1}(X)_\mathbb{R}\). Also, any current \(T_\omega\) associated to a real \((1,1)\)-form \(\omega\) on \(X\) is necessarily a \((1,1)\)-current; so every class in \(H^{1,1}(X)_\mathbb{R}\) is represented by some \(d\)-closed \((1,1)\)-current. A \((1,1)\)-current \(T\) is positive if it satisfies \(T(\eta) \geq 0\) for any semi-positive \(\eta \in \Omega^{1,1}(X)_\mathbb{R}\).
\begin{proposition}
Let \(X\) be a connected compact K\"{a}hler surface. Then the trivial current is the unique positive current on \(X\) representing the trivial class in \(H^{1,1}(X)_\mathbb{R}\).
\end{proposition}

\emph{Proof:} Suppose that \(T\) is a positive current representing the trivial class in \(H^{1,1}(X)_\mathbb{R}\). If \(T\) were non-trivial, then there would be some (non-closed) form \(\eta \in \Omega^{1,1}(X)_\mathbb{R}\) such that \(T(\eta) < 0\); but then there would be some \(\Kahler\) form \(\kappa \in \Omega^{1,1}(X)_\mathbb{R}\) such that
\[T(\eta + \kappa) = T(\eta) < 0\]
with \(\eta + \kappa\) a positive form,
which cannot happen. \(\Box\)
\\
\\
A \(4\)-current on \(X\) that is positive on non-negative functions is the same thing as a finite Borel measure on \(X\), and vice versa. The wedge product of a positive \((1,1)\)-current \(T\) with a semi-positive \((1,1)\)-form \(\omega\) is the measure \(T \wedge T_\omega\) given by, for a smooth function \(\phi \in \Omega^0(X)_\mathbb{R}\),
\[(T \wedge T_\omega)(\phi) = T(\phi \omega);\]
if \(T'\) is a positive \((1,1)\)-current that is a weak limit of semi-positive \((1,1)\)-forms, then the wedge product \(T \wedge T'\) is the weak limit of the wedge products of \(T\) with the semi-positive forms.
\begin{theorem}
[\cite{ca2} and \cite{ca1}] Let \(X\) be a connected compact K\"{a}hler surface, and suppose that \(\sigma\) is an automorphism of \(X\) with entropy \(\log(\lambda)>0\). Then there are positive closed \((1,1)\)-currents \(T_+\) and \(T_-\) on \(X\) with the following properties:
\begin{description}
\item[1)] \(\sigma^*T_+ = \lambda T_+\) and \(\sigma_*T_- = \lambda T_-\);
\item[2)] \(T_+\) and \(T_-\) are the unique positive \((1,1)\)-currents representing their respective cohomology classes;
\item[3)] If \(T\) is another positive \(d\)-closed \((1,1)\)-current on \(X\), then the sequences
\[\{\lambda^{-k}(\sigma^k)^*T\}_{k \in \mathbb{N}} \, \text{and} \, \{\lambda^{-k}(\sigma^k)_*T\}_{k \in \mathbb{N}}\]
converge weakly to \(([T_-].[T])T_+\) and \(([T_+].[T])T_-\), respectively; and
\item[4)] \(T_+ \wedge T_-\) is the unique measure of maximal entropy for \(\sigma\).
\end{description}
\end{theorem}
The currents \(T_+\) and \(T_-\) in Theorem 4.2 both have trivial wedge self-products. Thus, for any positive real numbers \(a\) and \(b\), the current \(aT_+ + bT_-\) has wedge self-product
\[a^2(T_+ \wedge T_+) + 2ab(T_+ \wedge T_-) + b^2(T_- \wedge T_-) = 2ab(T_+ \wedge T_-),\]
which is a positive multiple of the measure of maximal entropy for \(\sigma\). Moreover,
\[\sigma^*(aT_+ + bT_-) + (\sigma^{-1})^*(aT_+ + bT_-) = (\lambda + \lambda^{-1})(aT_+ + bT_-)\]
for any \(a\) and \(b\). The cohomology classes \([T_+]\) and \([T_-]\) are eigenvectors for the action of \(\sigma^*\) on \(H^{1,1}(X)_\mathbb{R}\) and are contained in \(\overline{C_K(X)}\); also, \([T_+].[T_-]=1\). The measures of maximal entropy for \(\sigma\), \(\sigma^{-1}\), and all iterates of \(\sigma\) or \(\sigma^{-1}\) are the same.
\begin{proposition}
Let \(X\) be a smooth complex projective surface, and suppose that \(\sigma\) is an automorphism of \(X\) with entropy \(\log(\lambda)>0\). Let \(L\) be a nef and big line bundle on \(X\), and suppose that \(\xch(L)\) contains some semi-positive form \(\omega_0\). Then the inductively defined sequence
\[\{\omega_n = (\lambda + \lambda^{-1})^{-1} (\sigma^*\omega_{n-1} + (\sigma^{-1})^*\omega_{n-1})\}_{n \in \mathbb{N}}\]
converges weakly to a current \(T\) with the property that \(T \wedge T\) is some positive scaling of the measure of maximal entropy for \(\sigma\). 
\end{proposition}

\emph{Proof:} Since \(\omega_0\) is semi-positive and \(d\)-closed, the current \(T_{\omega_0}\) is positive and \(d\)-closed; also, \(\xch(L)\) and \([T_{\omega_0}]\) are the same cohomology class. Since \([T_+]\) and \([T_-]\) are contained in \(\overline{C_K(X)} \subseteq \overline{C_+(X)}\), it follows from the analogue of the Hodge index theorem that the Chern class of any nef and big line bundle on \(X\) must have positive intersection with both \([T_+]\) and \([T_-]\). So, in particular, the sequence \[\{T_n = \lambda^{-n}((\sigma^n)^*T_{\omega_0}+((\sigma^{-1})^n)^*T_{\omega_0})\}_{n \in \mathbb{N}}\]
converges weakly to \(aT_+ + bT_-\) for some positive real numbers \(a\) and \(b\). For any \(k \in \mathbb{N}\),
\[T_{\omega_{2k}} = \left( \frac{1}{\lambda + \lambda^{-1}} \right)^{2k} \sum\nolimits_{j=0}^k
\left( \begin{array}{c}
2k \\
j \end{array} \right)
\lambda^{2k-2j} T_{2k-2j} \]
(where \(T_0 = T_{\omega_0}\)) and
\[T_{\omega_{2k-1}} = \left( \frac{1}{\lambda + \lambda^{-1}} \right)^{2k-1} \sum\nolimits_{j=0}^{k-1}
\left( \begin{array}{c}
2k-1 \\
j \end{array} \right)
\lambda^{2k-1-2j} T_{2k-1-2j};\]
also, for any \(m \in \mathbb{N}\),
\[\sum\nolimits_{j=0}^m
\left( \begin{array}{c}
m \\
j \end{array} \right)
\lambda^{m-2j} = (\lambda + \lambda^{-1})^m > 2^m = \sum\nolimits_{j=0}^m
\left( \begin{array}{c}
m \\
j \end{array} \right) .\]
(Compare \cite{ka1}, \S 3.3.2.) Moreover,
\[\lim_{k \rightarrow \infty} \left( \frac{1}{\lambda + \lambda^{-1}} \right)^{2k}
\left( \begin{array}{c}
2k \\
k-l \end{array} \right)
\lambda^{2l} = \lim_{k \rightarrow \infty} \left( \frac{1}{\lambda + \lambda^{-1}} \right)^{2k-1}
\left( \begin{array}{c}
2k-1 \\
k-l \end{array} \right)
\lambda^{2l-1} = 0\]
for any \(l \in \mathbb{N}_0\), and
\[\lim_{m \rightarrow \infty} \left( \frac{1}{\lambda + \lambda^{-1}} \right)^m \sum\nolimits_{m/2 < j \leq m}
\left( \begin{array}{c}
m \\
j \end{array} \right)
\lambda^{m-2j} = 0.\]
Thus it follows that, for any real \((1,1)\)-form \(\eta\), the sequences
\(\{T_{\omega_{2k}}(\eta)\}_{k \in \mathbb{N}}\) and \(\{T_{\omega_{2k-1}}(\eta)\}_{k \in \mathbb{N}}\)
both converge to \(aT_+(\eta) + bT_-(\eta)\). \(\Box\)

\subsection{Periodic Subvarieties}

It is a well-known fact that the entropy of any automorphism of a connected compact K\"{a}hler curve is zero. (See, e.g., \cite{ca1}, \S 1.1 and \S 4.4.) Suppose that \(X\) is projective, and that \(\sigma\) has entropy \(\log(\lambda) > 0\). Let \(\mu_\sigma\) be the measure of maximal entropy for \(\sigma\). By definition, the entropy of \(\sigma\) restricted to the orbit of any periodic point must be zero. Suppose that \(C\) is an irreducible curve on \(X\) that is periodic for \(\sigma\); so \(C\) is a fixed curve for \(\sigma^k\) for some \(k \in \mathbb{N}\). A finite sequence of blow-ups of singular points of \(C\) yields a surface \(X'\) in which \(C'\), the strict transform of \(C\), is non-singular. (See \cite{bhpv}, \S II.7.) Since \(\sigma^k\) must preserve the set of non-singular points on \(C\), each blown-up point in the construction of \(X'\) must be periodic for \(\sigma^k\); so \(\sigma^k\) extends to an automorphism of \(X'\), which in turn restricts to an automorphism (with entropy zero) of the non-singular curve \(C'\). Thus the restriction of \(\sigma^k\) to \(C\) must have entropy zero. Since the entropy of any iterate of \(\sigma\) restricted to some \(\sigma\)-invariant set is at least the entropy of \(\sigma\) on the set, the entropy of \(\sigma\) restricted to the orbit of any periodic curve must be zero.

\begin{proposition}
Let \(X\) be a smooth complex projective surface, and let \(\sigma\) be an automorphism of \(X\). Suppose that \(\sigma\) has positive entropy, and that \(C \subseteq X\) is a proper irreducible subvariety (i.e., a point or an irreducible curve) that is periodic for \(\sigma\). Then \(\mu_\sigma(C) = 0\).
\end{proposition}

\emph{Proof:} Let \(k \in \mathbb{N}\) be the period of \(C\) under \(\sigma\); so the entropy of \(\sigma^k\) restricted to \(C\) is zero. Let \(\nu\) be the measure on \(X\) given by, for a Borel set \(A\),
\[\nu(A) = (1-\mu_\sigma(C))^{-1}\mu_\sigma(A-(A \cap C)).\]
If \(\mu_\sigma(C)\) were positive, then \(\nu\) would be a \(\sigma^k\)-invariant probability measure on \(X\) such that the entropy of \(\sigma^k\) with respect to \(\nu\) is strictly greater than the entropy of \(\sigma^k\) with respect to \(\mu_\sigma\)--which cannot exist. \(\Box\)
\\
\\
If the support of \(\mu_\sigma\) were contained in some proper subvariety of \(X\), then it would be contained in some finite union of points and irreducible curves on \(X\), each periodic for \(\sigma\); but then the support of \(\mu_\sigma\) would have measure zero, which cannot be the case. Thus the support of \(\mu_\sigma\) is Zariski dense in \(X\). If \(\mu_\sigma(\{x\})\) were positive for some point \(x \in X\), then \(x\) would necessarily be a periodic point--and hence could not have positive measure; so \(\mu_\sigma\) has no atoms on \(X\). The same argument shows that the measure of any irreducible curve on \(X\) must in fact be zero.
\\
\\
For any \(k \in \mathbb{N}\) , the set of points in \(X\) of exact period \(k\) for \(\sigma\) is a subvariety \(Z_k \subseteq X\); the set \(\xPer(\sigma,k)\) of isolated points of exact period \(k\) is the complement in \(Z_k\) of the union of all of the curves contained in \(Z_k\). Since there can only be finitely many curves on \(X\) that are periodic for \(\sigma\), there is an upper bound on the values of \(k\) for which \(\xPer(\sigma,k)\) is a proper subset of \(Z_k\).
\begin{theorem}
[\cite{ca2} and \cite{ca1}] Let \(X\) be a smooth complex projective surface, and suppose that \(\sigma\) is an automorphism of \(X\) with positive entropy. Then the isolated periodic points for \(\sigma\) are equidistributed with respect to the measure of maximal entropy for \(\sigma\), in the sense that the sequence
\[\left\{ (|\xPer(\sigma,n)|)^{-1} \sum\nolimits_{x \in \xPer(\sigma,n)} \delta_x \right\}_{n \in \mathbb{N}}\]
converges weakly to \(\mu_\sigma\).
\end{theorem}
Since the support of \(\mu_\sigma\) is Zariski dense in \(X\), the set of (isolated) points in \(X\) that are periodic for \(\sigma\) must also be Zariski dense. (See also \cite{ca1}, \S 4.4.3.)

\subsection{Semi-Positive Forms}

If a line bundle \(L\) on a smooth complex projective surface has a semi-positive form \(\omega\) in its Chern class, then the intersection properties of \(\omega\) guarantee that \(L\) is nef. The following theorem (primarily due to Kawamata) provides a means of finding line bundles whose Chern classes contain semi-positive forms. 
\begin{theorem}
[\cite{ein}, Theorem 2, and \cite{huy}, Corollary 4.3.19] Let \(A\) be a nef and big divisor on a smooth complex projective surface \(X\), and suppose that \(K_X \otimes A\) is nef. Then \((K_X \otimes A)^{\otimes q}\) has a semi-positive form in its Chern class for some \(q \in \mathbb{N}\).
\end{theorem}
The line bundle \((K_X \otimes A)^{\otimes q}\) in Theorem 4.6 is in fact a base-point-free effective divisor class, and it is this feature that guarantees the existence of a semi-positive form in its Chern class.
\begin{proposition}
Let \(X\) be a smooth complex projective surface, let \(\sigma\) be an automorphism of \(X\) with positive entropy, and let \(L\) be a nef and big distinguished line bundle on \(X\). Suppose that \(X\) is birational to a K3 surface, an Enriques surface, or an abelian surface. Then there is some \(q \in \mathbb{N}\) such that \(L^{\otimes q}\) has a semi-positive form in its Chern class.
\end{proposition}

\emph{Proof:} If \(X\) is a K3 surface or an abelian surface, then \(K_X\) is trivial; so
\[K_X \otimes L = L\]
is nef (and big), and \(\xch(L^{\otimes q})\) contains a semi-positive form for some \(q \in \mathbb{N}\). (See \cite{bhpv}, \S VI.1.) If \(X\) is an Enriques surface, then \(K_X \otimes K_X\) is trivial and \(K_X\) must have zero intersection with every line bundle on \(X\); so \(K_X \otimes L\) is nef (and big), and
\[(K_X \otimes L)^{\otimes q} = K_X^{\otimes (q \bmod{2})} \otimes L^{\otimes q}\]
has a semi-positive form in its Chern class for some \(q \in \mathbb{N}\)--from which it follows that \(\xch(L^{\otimes 2q})\) contains a semi-positive form. (See \cite{bhpv}, \S VIII.15.)

More generally (if \(X\) is not minimal), let \(\mathcal{F}\) be the set of \((-1)\)-curves on \(X\); so \(\mathcal{F}\) is a finite set, and every \(E \in \mathcal{F}\) is periodic for \(\sigma\). Let \(X'\) be the surface obtained from \(X\) by contraction of all curves in \(\mathcal{F}\); so the Picard group of \(X\) is given by
\[\xPic(X) \cong \xPic(X') \times \left( \bigoplus\nolimits_{E \in \mathcal{F}} \mathopen{<}[E]\mathclose{>} \right).\]
Since every \(E \in \mathcal{F}\) must have zero intersection with \(L\), \(L\) can be expressed as
\[L = (L',0) \in \xPic(X),\]
for some nef and big \(L' \in \xPic(X')\). Since \(X'\) is a K3 surface, an Enriques surface, or an abelian surface, \((L')^{\otimes q}\) is a base-point-free effective divisor class for some \(q \in \mathbb{N}\). For any \(x \in X\), let \(x'\) be the image of \(x\) in \(X'\); then there is some effective divisor \(D_x\) representing \(L'\) that does not contain \(x'\). Thus, for any \(x \in X\), the proper transform of \(D_x\) in \(X\) is an effective divisor representing \(L\) that does not contain \(x\); so \(L^{\otimes q}\) is a base-point-free effective divisor class and hence has a semi-positive form in its Chern class. \(\Box\)
\\
\\
\emph{Proof of Theorem 1.3:} By Proposition 3.6, there is a nef and big distinguished line bundle \(L \in \xPic(X)\); if \(X\) is not rational, then, by Proposition 4.7, some multiple \(L^{\otimes q}\) is a distinguished nef and big line bundle whose Chern class contains a semi-positive form; if \(\sigma\) has no periodic curves, then, by Proposition 3.5, \(L\) can be taken to be ample so that its Chern class contains a positive form. The weak convergence of the sequence and the properties of the limit are given by Proposition 4.3. \(\Box\)

\section{Orbits of Line Bundles}

In this section, we investigate the orbits of Chern classes under projective surface automorphisms with positive entropy. We prove one direction of Theorem 1.4 in \S 5.1, and we complete the proof in \S 5.2.
\\
\\
Let \(X\) be a smooth complex projective surface, and let \(\sigma\) be an automorphism of \(X\) with entropy \(\log(\lambda) > 0 \). Let \(S(t)\) be the minimal polynomial for \(\lambda\), and let \(s\) be the degree of \(S(t)\). For any line bundle \(L \in \xPic(X)\), the span of the orbit of \(\xch(L)\) under \(\sigma^*\) is a subspace of \(\xNS(X)_\mathbb{R} = \xNS(X) \otimes \mathbb{R}\) that is invariant under \(\sigma^*\). Let \(\xNS'(X)\) be the \(s\)-dimensional sublattice of \(\xNS(X)\) consisting of all Chern classes that are annihilated by \(S(\sigma^*)\); so the characteristic polynomial for the action of \(\sigma^*\) on \(\xNS'(X)_\mathbb{R} = \xNS'(X) \otimes \mathbb{R}\) is \(S(t)\). Let \(C(t)\) be the characteristic polynomial for the action of \(\sigma^*\) on \(\xNS'(X)_\mathbb{R}^\perp \subset \xNS(X)_\mathbb{R}\); so \(C(t)\) is either trivial (if \(\xNS'(X) = \xNS(X)\)) or a product of cyclotomic polynomials.

\subsection{Inseparable Cyclotomic Factors}

It is a well-known fact that any orthogonal operator (on a finite-dimensional space) is diagonalizable over \(\mathbb{C}\).
\begin{comment}
Then the quadratic form extends to a definite hermitian form on \(W \otimes \mathbb{C}\), and \(T\) extends to a linear transformation \(T_\mathbb{C}\) of \(W \otimes \mathbb{C}\) that preserves the hermitian form. Let \(\{w_1, \dots , w_r\}\) be a basis for \(W \otimes \mathbb{C}\) that expresses \(T_\mathbb{C}\) in Jordan canonical form. If \(T_\mathbb{C}\) is not diagonal with respect to this basis, then there is some basis element \(w_j\) such that \(T_\mathbb{C}w_j = w_{j-1} + \eta w_j\) and \(T_\mathbb{C}w_{j-1} = \eta w_{j-1}\), where \(\eta\) is some non-simple eigenvalue of \(T\). Since \((T_\mathbb{C}w_{j-1})^2 = |\eta|^2(w_{j-1})^2\) and \((w_{j-1})^2 \neq 0\), \(\eta\) must lie on the unit circle; but then
\[(T_\mathbb{C}w_{j-1}).(T_\mathbb{C}w_j) = (\eta w_{j-1}).(w_{j-1}+\eta w_j) = \eta (w_{j-1})^2 + w_{j-1}.w_j \]
with \(\eta(w_{j-1})^2 \neq 0\), which cannot be the case. So T must be diagonalizable over \(\mathbb{C}\).
\end{comment}
\begin{proposition}
Assume the hypotheses and notation of the preceding text in this section. Then, as a linear transformation of \(\xNS(X)_\mathbb{R}\), \(\sigma^*\) is diagonalizable over \(\mathbb{C}\).
\end{proposition}

\emph{Proof:} Since \(S(t)\) is separable, each of its roots is a simple eigenvalue of \(\sigma^*\); thus \(\sigma^*\) is diagonalizable over \(\mathbb{C}\) as a linear transformation of \(\xNS'(X)_\mathbb{R}\). Since \(\xNS'(X)_\mathbb{R}^\perp \subseteq \xNS(X)_\mathbb{R}\) is either negative definite (so that \(\sigma^*\) is an orthogonal operator on \(\xNS'(X)_\mathbb{R}^\perp\)) or trivial, \(\sigma^*\) is in fact diagonalizable over \(\mathbb{C}\) as a linear transformation of \(\xNS(X)_\mathbb{R}\). \(\Box\)
\\
\\
If a linear transformation of a complex vector space is diagonalizable, then its minimal polynomial is necessarily separable.
\begin{proposition}
Assume the hypotheses and notation of the preceding text in this section, and suppose that \(C(t)\) is inseparable. Then no line bundle on \(X\) can have an orbit of Chern classes that spans \(\xNS(X)_\mathbb{R}\).
\end{proposition}

\emph{Proof:} Let \(M(t) \in \mathbb{Z}[t]\) be the minimal polynomial for the action of \(\sigma^*\) on \(\xNS(X)_\mathbb{R}\), and let \(m \in \mathbb{N}\) be the degree of \(M(t)\). So the orbit of any Chern class \(\xch(L) \in \xNS(X)\) must be contained in
\[\mathopen{<} \xch(L),\dots,(\sigma^*)^{m-1}\xch(L) \mathclose{>}.\]
Since \(C(t)\) is inseparable, there is a cyclotomic polynomial \(P(t)\) such that \(P(t)^2\) divides \(C(t)\); however, \(P(t)^2\) cannot divide \(M(t)\). It follows that \(m\) is strictly less than the degree of \(S(t)C(t)\); so the orbit of any Chern class in \(\xNS(X)\) must be contained in
a subspace of dimension strictly less than that of \(\xNS(X)_\mathbb{R}\). \(\Box\)
\\
\\
If the degree of \(C(t)\) is two or more, then \((x-1)^2\) will necessarily divide the cyclotomic factor for some iterate of \(\sigma\) (so this cyclotomic factor will be inseparable). Also, the degree of \(C(t)\) must be at least the number of curves on \(X\) that are periodic for \(\sigma\).

\subsection{Separable Cyclotomic Factors}

Since \(X\) is projective, the intersection \(\xNS(X) \cap C_K(X)\) is non-empty; moreover, since \(C_K(X)\) is open in \(H^{1,1}(X)_\mathbb{R}\), it is also open in \(\xNS(X)_\mathbb{R}\). It follows that \(\xNS(X)_\mathbb{R}\) is spanned by the set of all Chern classes of ample line bundles on \(X\).
\\
\\
\emph{Proof of Theorem 1.4:} Let \(Q(t) \in \mathbb{Z}[t]\) be the characteristic polynomial for the action of \(\sigma^*\) on \(\xNS(X)\). If \(Q(t)\) is inseparable, then, by Proposition 5.2, no line bundle in \(\xPic(X)\) has an orbit of Chern classes that spans \(\xNS(X) \otimes \mathbb{R}\). If \(Q(t)\) is separable, then it is also the minimal polynomial for the action of \(\sigma^*\) on \(\xNS(X)\); in this case, the set of elementary divisors \(\{Q_1(t),\dots,Q_r(t)\}\) for \(\sigma^*\) is precisely the set of irreducible factors of \(Q(t)\), and Lemma 3.2.1 shows that there is a decomposition of \(\xNS(X) \otimes \mathbb{Q}\) into spaces \(E_1, \dots, E_r\) such that each \(E_j\) is invariant under \(\sigma^*\) and annihilated by \(Q_j(\sigma^*)\). Let \(\{A_1,\dots,A_q\}\) be a set of ample line bundles whose Chern classes span \(\xNS(X) \otimes \mathbb{R}\); then there is some linear combination
\[A = k_1A_1 + \dots + k_qA_q\]
with coefficients in \(\mathbb{N}\) such that the projection of \(\xch(A)\) to any \(E_j\) is not zero. So \(A\) is an ample line bundle whose orbit of Chern classes spans \(\xNS(X) \otimes \mathbb{R}\). \(\Box\)
\\
\\
If \(C(t)\) is trivial or linear, then every iterate of \(\sigma\) will have an induced action on \(\xNS(X)\) that is separable.

\section{Examples of Surface Automorphisms}

In each of the examples that follow, the surface is projective and the automorphism is minimal in the sense that no exceptional curve (of the first kind) is periodic. The examples highlight:
\begin{description}
\item[1)] Automorphisms with and without periodic curves;
\item[2)] Automorphisms with and without separable cyclotomic factors; and
\item[3)] Automorphisms with quadratic entropies.
\end{description}

\subsection{Abelian Surfaces and Kummer Surfaces}

An abelian surface automorphism with positive entropy cannot have a periodic curve (since any such curve would be rational). (See \cite{bhpv}, \S II.11, and \cite{ka2}, \S 3.) The possible positive values of the entropy of an abelian surface automorphism are limited to logarithms of Salem numbers of degree two or four.
\\
\\
Any automorphism of a two-dimensional complex torus descends to an automorphism of the associated Kummer surface; the sixteen curves on the Kummer surface arising from the fixed points for the involution on the torus must be periodic for the Kummer surface automorphism, and the entropy of the Kummer surface automorphism must be the same as the entropy of the torus automorphism. (See also \cite{mc1}, \S 4.) A Kummer surface may admit additional automorphisms beyond those which come from automorphisms of the torus to which the Kummer surface is associated. (See, e.g., \cite{keu}.)
\\
\\
\emph{Abelian surface automorphisms with quadratic entropies.} For any elliptic curve \(E = \mathbb{C}/\Lambda\), the product \(E \times E = \mathbb{C}^2/(\Lambda \times \Lambda)\) is an abelian surface. Via the group law on \(E\), any element \(A = (a_{ij}) \in \xGL_2(\mathbb{Z})\) gives an automorphism of \(E \times E\) by, for \((e_1,e_2) \in E \times E\),
\[A(e_1,e_2)=(a_{11}e_1+a_{12}e_2,a_{21}e_1+a_{22}e_2).\]
The entropy of \(A\) is the maximum of the squares of the magnitudes of its eigenvalues; so the entropy is the logarithm of a degree-two Salem number whenever the eigenvalues are not both units. (See also \cite{res}, \S 5.) If \(E\) is chosen generically, then the Picard rank of \(E \times E\) is three, and any automorphism of \(E \times E\) with positive entropy must have as its entropy the logarithm of a degree-two Salem number and hence also must have a separable cyclotomic factor.
\\
\\
\emph{Abelian surface automorphisms with degree-four entropies.} If \(E\) in the previous example admits multiplication by some \(\omega \in \mathbb{C}-\mathbb{R}\), then \(E \times E\) has Picard rank four. Here, elements of \(\xGL_2(\mathbb{Z}[\omega])\) give automorphisms of \(E \times E\), and the entropies of some of these automorphisms are logarithms of degree-four Salem numbers (in which case the automorphisms have trivial cyclotomic factors).

\subsection{Non-Kummer K3 Surfaces}

\emph{Automorphisms of K3 surfaces that do not contain any smooth rational curves.} Any smooth hypersurface in \(\mathbb{P}^1 \times \mathbb{P}^1 \times \mathbb{P}^1\) that is given by an equation of tri-degree \((2,2,2)\) is a K3 surface; any such surface admits three non-commuting involutions whose compositions give automorphisms with positive entropy. (See \cite{maz}, \S 9, and \cite{ca2}, \S 1.5.) If the defining equation is chosen generically, then the Picard rank of the K3 surface is three, and the intersection form on the Picard group is given as a matrix by
\[\left( \begin{array}{ccc}
0 & 2 & 2 \\
2 & 0 & 2 \\
2 & 2 & 0 \end{array} \right).\]
(See also \cite{bar}, \S 2.) So no line bundle on a generic K3 surface in this family can have self-intersection equal to \(-2\), and thus no curve on a generic K3 surface in this family can have negative self-intersection; so any positive-entropy automorphism of such a surface has no periodic curves. Moreover, any such automorphism must have as its entropy the logarithm of a degree-two Salem number, and also must have a separable cyclotomic factor. (See also \cite{lee}, \S 2.)
\\
\\
\emph{K3 surface automorphisms with infinite orbits of smooth rational curves but no periodic curves.} If the equation in the previous example is, in affine coordinates,
\[f(x,y,z) = (x^2+1)(y^2+1)(z^2+1) + Axyz - 2,\]
for some \(A \in \mathbb{C}\), then the K3 surface has twelve (smooth rational) curves with self-intersection \(-2\). (See \cite{mc1}, \S 1, and \cite{row}.) The positive-entropy automorphisms from the previous example still have no periodic curves (so the orbits of the smooth rational curves are infinite).

If the equation in the previous example is, in affine coordinates,
\[f(x,y,z) = x^2(y^2 + yz +z^2 + z) + x(y^2z^2 + y^2z + z) + (y^2z^2 +y^2z + y + z),\]
then the K3 surface contains the smooth rational curve \(\mathbb{P}^1 \times \{0\} \times \{0\}\), and the Picard rank of the K3 surface is four; moreover, there is an additional involution on the surface that does not commute with the three involutions from the previous example. (See \cite{bar}, \S 2.) For all  compositions of the four involutions that have been investigated, the positive-entropy automorphisms have no periodic curves (so the orbit of the projective line is infinite) and have entropies that are logarithms of degree-two Salem numbers. (See also \cite{lee}, \S 2.) These examples include some automorphisms with separable but not irreducible degree-two cyclotomic factors.
\\
\\
\emph{A non-Kummer K3 surface automorphism with periodic curves.} The Torelli theorem leads to synthetic examples of K3 surface automorphisms with many interesting properties. In particular, there is a projective K3 surface automorphism whose entropy is the logarithm of a Salem number of degree six, and whose Picard group contains a two-dimensional sublattice on which the automorphism induces the identity map and the intersection form is given as a matrix by
\[\left( \begin{array}{cc}
-2 & -1 \\
-1 & -20 \end{array} \right).\]
(See \cite{mc4}, \S 8.) The Riemann-Roch theorem implies that there is a periodic rational curve for this automorphism. (See \cite{bar}, \S 1.)
\\
\\
\emph{K3 surface automorphisms with entropies of degree equal to Picard rank.} Another synthetic construction gives a projective K3 surface automorphism in which the degree of the entropy and the Picard rank of the surface are both eighteen (so the automorphism has a trivial cyclotomic factor). (See \cite{mc3}, \S 8.) In this case, the automorphism cannot have a periodic curve, since no root of unity is can be eigenvalue for the action on the N\'eron-Severi group induced by the automorphism.

Similarly, there are projective K3 surfaces with Picard rank two that admit automorphisms with positive entropy. (See \cite{ogu}, \S 4.) The automorphisms here must have quadratic entropies, and also must have trivial cyclotomic factors.

\subsection{Rational Surfaces}

\emph{Rational surface automorphisms with no periodic curves.} For any two complex numbers \(\alpha\) and \(\beta\), the map \(f_{\alpha,\beta}\) given by
\[f_{\alpha,\beta}(x,y) = (y, (y+\alpha)/(x+\beta))\]
defines a birational self-map of \(\mathbb{P}^2\). For certain parameters, a finite blow-up resolves the map \(f_{\alpha,\beta}\) into a rational surface automorphism with positive entropy. Some of the automorphisms with positive entropy constructed in this way have periodic curves and some do not. (See \cite{be2}.)
\\
\\
\emph{Rational surface automorphisms with invariant curves.} Most of the known examples of rational surface automorphisms with positive entropy arise from blow-ups of \(\mathbb{P}^2\) at points along invariant curves. (See \cite{be1}, \cite{dil}, \cite{mc2}, and \cite{ueh}.) Indeed, any positive value that is the spectral radius of an element in the Weyl group for some Lorentz-Minkowski lattice is also the exponent of the entropy of some automorphism of a blow-up of \(\mathbb{P}^2\) with an invariant anticanonical curve. (See \cite{ueh}, \S 1.) In general, the constructions of these rational surface automorphisms yield entropies that are logarithms of Salem numbers of arbitrarily large degrees, as well as arbitrarily large cyclotomic factors; however, it can happen that such an automorphism has a linear cyclotomic factor. (See \cite{mc2}, \S 1 and \S 2.)
\\
\\
\noindent \emph{A rational surface automorphism with quadratic entropy.} The matrix
\[\left( \begin{array}{ccccccccccccc}
92 & 8 & 8 & 8 & 14 & 14 & 14 & 25 & 25 & 25 & 44 & 44 & 44 \\
-8 & 0 & -1 & -1 & -1 & -1 & -1 & -2 & -2 & -2 & -4 & -4 & -4 \\
-8 & -1 & 0 & -1 & -1 & -1 & -1 & -2 & -2 & -2 & -4 & -4 & -4 \\
-8 & -1 & -1 & 0 & -1 & -1 & -1 & -2 & -2 & -2 & -4 & -4 & -4 \\
-44 & -4 & -4 & -4 & -6 & -7 & -7 & -12 & -12 & -12 & -21 & -21 & -21 \\
-44 & -4 & -4 & -4 & -7 & -6 & -7 & -12 & -12 & -12 & -21 & -21 & -21 \\
-44 & -4 & -4 & -4 & -7 & -7 & -6 & -12 & -12 & -12 & -21 & -21 & -21 \\
-25 & -2 & -2 & -2 & -4 & -4 & -4 & -6 & -7 & -7 & -12 & -12 & -12 \\
-25 & -2 & -2 & -2 & -4 & -4 & -4 & -7 & -6 & -7 & -12 & -12 & -12 \\
-25 & -2 & -2 & -2 & -4 & -4 & -4 & -7 & -7 & -6 & -12 & -12 & -12 \\
-14 & -1 & -1 & -1 & -2 & -2 & -2 & -4 & -4 & -4 & -6 & -7 & -7 \\
-14 & -1 & -1 & -1 & -2 & -2 & -2 & -4 & -4 & -4 & -7 & -6 & -7 \\
-14 & -1 & -1 & -1 & -2 & -2 & -2 & -4 & -4 & -4 & -7 & -7 & -6 \end{array} \right) \]
gives a Weyl element for the Lorentz-Minkowski lattice of dimension thirteen whose spectral radius is a degree-two Salem number; thus there is an automorphism of \(\mathbb{P}^2\) blown up at twelve points whose entropy is the logarithm of this degree-two Salem number. (This explicit example was constructed via trial and error as a composition of known Weyl elements, based in part on ideas in \cite{ueh}.)

\bibliographystyle{plain}
\bibliography{ReschkeP-refs-2014.01.30}

\begin{comment}

\end{comment}

\end{document}